\documentclass[10.5pt]{amsart}
\usepackage{graphicx}
\usepackage{amssymb}
\usepackage{amsmath}
\usepackage{amsthm}
\usepackage{amscd}
\usepackage[all,2cell,dvips]{xy} \UseAllTwocells
\SilentMatrices\newtheorem{thm}{Theorem}[section]

\newtheorem{lem}[thm]{Lemma}

\newtheorem{prop}[thm]{Proposition}
\theoremstyle{definition}

\theoremstyle{remark}

\numberwithin{equation}{section}

\begin{document}
\title[Algebras of derived dimension 0]{Algebras of derived dimension zero}
\author[Xiao-Wu Chen, Yu Ye \ and\   Pu Zhang
] {Xiao-Wu Chen$^A$, Yu Ye$^A$ \ and \  Pu Zhang$^{B,*}$}
\thanks{$^*$ The corresponding author}
\thanks{Supported in part by the National Natural Science Foundation of China (Grant No. 10301033)}
\thanks{Email: xwchen$\symbol{64}$mail.ustc.edu.cn, yeyu$\symbol{64}$ustc.edu.cn, pzhang$\symbol{64}$sjtu.edu.cn}
%
\maketitle
\begin{center}
$^A$Department of Mathematics \\ University of Science and
Technology of China \\Hefei 230026, Anhui, P. R. China
\end{center}
\begin{center}
$^B$Department of Mathematics, \ \ Shanghai Jiao Tong University\\
Shanghai 200240, P. R. China
\end{center}
\begin{center}
\end{center}
\section{Introduction}
\subsection{} A dimension for a triangulated
category has been introduced by Rouquier in [Ro], which gives a new
invariant for algebras and algebraic varieties under derived
equivalences. For related topics see also [BV], and [Hap],
p.70.\vskip5pt

Let $\mathcal{C}$ be a triangulated category with shift functor
$[1]$,  $\mathcal{I}$ and $\mathcal{J}$ full subcategories of
$\mathcal{C}$. Denote by $\langle \mathcal{I} \rangle$ the smallest
full subcategory of $\mathcal{C}$ containing $\mathcal{I}$ and
closed under isomorphisms, finite direct sums, direct summands, and
shifts. Any object of $\langle \mathcal{I}\rangle$ is isomorphic to
a direct summand of a finite direct sum $\bigoplus\limits_i
I_i[n_i]$ with each $I_i \in \mathcal{I}$ and $n_i\in\mathbb{Z}$.
Define $\mathcal{I} \star \mathcal{J}$ to be the full subcategory of
$\mathcal{C}$ consisting of the objects $M$, for which there is a
distinguished triangle $I \longrightarrow M \longrightarrow J
\longrightarrow I[1]$ with $I \in \mathcal{I}$ and $J \in
\mathcal{J}$. Now define $\langle \mathcal{I} \rangle_0: =\{0\}$,
 and
$\langle \mathcal{I}\rangle_n: =\langle \langle
\mathcal{I}\rangle_{n-1} \star \langle \mathcal{I}\rangle \rangle
$ for $n \geq 1$. Then $\langle \mathcal{I} \rangle_1 = \langle
\mathcal{I} \rangle $, and  $\langle \mathcal{I} \rangle_n=\langle
\langle \mathcal{I} \rangle \star \cdots \star \langle \mathcal{I}
\rangle \rangle$, by the associativity of $\star$ (see [BV]). Note
that $\langle \mathcal{I} \rangle_{\infty}: =\bigcup_{n=0}^\infty
\langle \mathcal{I} \rangle_n$ is the smallest thick triangulated
subcategory of $\mathcal{C}$ containing $\mathcal{I}$.\par \vskip
5pt

\vskip5pt

By definition, the \textit{dimension} of $\mathcal{C}$, denoted by
${\rm dim}(\mathcal {C})$,  is the minimal integer $d\ge 0$ such
that there exists an object $M \in \mathcal{C}$ with
$\mathcal{C}=\langle M \rangle_{d+1}$, or $\infty$ when there is no
such an object $M$. See [Ro].

\vskip5pt

Let  $A$ be a finite-dimensional algebra over a field $k$. Denote
by $A\mbox{-mod}$ the category of finite-dimensional left
 $A$-modules, and by $D^b(A\mbox{-mod})$ the bounded derived
category. Define the \textit{derived dimension} of $A$, denoted by
${\rm der.dim}(A)$, to be the dimension of the triangulated
category $D^b(A\mbox{-mod}).$ By [Ro] and [KK] one has
$${\rm der.dim}(A)\leq {\rm min}\{l(A), \; {\rm gl.dim}(A), \; {\rm
rep.dim}(A)\}$$
 where
$l(A)$ is the smallest integer $l\ge 0$ such that ${\rm
rad}^{l+1}(A)=0$, ${\rm gl.dim}(A)$ and ${\rm rep.dim}(A)$ are the
global dimension and the representation dimension of $A$ (for the
definition of ${\rm rep.dim}(A)$ see [Au]), respectively. In
particular, we have ${\rm der.dim}(A)< \infty$.\par \vskip 5pt

Our main result is \vskip 5pt

 \noindent{\bf Theorem}\; \emph{Let $A$ be a finite-dimensional algebra
over an algebraically closed field $k$. Then ${\rm der.dim}(A)=0$ if
and only if $A$ is an iterated tilted algebra of Dynkin type.}

\vskip 10pt

\subsection{} Let us fix some notation. For an
additive category $\mathcal{A}$, denote by $C^*(\mathcal{A})$ the
category of complexes of $\mathcal{A}$, where $*\in \{-, +, b\}$
means bounded-above, bounded-below, and bounded, respectively; and
by $C(\mathcal{A})$ the category of unbounded complexes. Denote by
$K^*(\mathcal{A})$ the corresponding homotopy category. If
$\mathcal{A}$ is abelian, we have derived category
$D^*(\mathcal{A})$.\par

For a finite-dimensional algebra $A$, denote by $A\mbox{-mod}$,
$A\mbox{-proj}$ and $A\mbox{-inj}$ the category of
finite-dimensional left $A$-modules, projective $A$-modules and
injective $A$-modules, respectively.\par

For triangulated categories and derived categories we refer to [V],
[Har], and [Hap]; for representation theory of algebras we refer to
[ARS] and [Ri]; and for tilting theory we refer to [Ri] and [Hap],
in particular, for iterated tilted algebras we refer to [Hap],
p.171.

\section{Proof of Theorem}

Before giving the proof of Theorem, we make some preparations.

\subsection{}  Let $A=\bigoplus_{j\geq 0} A_{(j)}$ be a
finite-dimensional positively-graded algebra over $k$, and
$A\mbox{-gr}$ the category of finite-dimensional left
$\mathbb{Z}$-graded $A$-modules with morphisms of degree zero. An
object in $A\mbox{-gr}$ is written as $M=\bigoplus_{j\in \mathbb{Z}}
M_{(j)}$. For each $i\in \mathbb{Z}$, we have the degree-shift
functor $(i): A\mbox{-gr} \longrightarrow A\mbox{-gr}$, defined by
$M(i)_{(j)}=M_{(i+j)}$, $\forall \ j \in \mathbb{Z}$. Let $U:
A\mbox{-gr} \longrightarrow A\mbox{-mod}$ be the degree-forgetful
functor. Then $U(M(i))=U(M)$, $\forall \ i\in \mathbb{Z}$. Clearly,
$A\mbox{-gr}$ is a Hom-finite abelian category, and hence by Remark
A.2 in Appendix it is Krull-Schmidt. An indecomposable in
$A\mbox{-gr}$ is called a gr-indecomposable module. The category
$A\mbox{-gr}$ has projective covers and injective hulls. Assume that
$\{e_1, e_2, \cdots, e_n\}$ is a set of orthogonal primitive
idempotents of $A_{(0)}$, such that $\{P_i:=Ae_i=\bigoplus_{j\geq 0}
A_{(j)}e_i\ | \ 1\le i\le n\}$  is a complete set of pairwise
non-isomorphic indecomposable projective $A$-modules. Then $P_i$
(resp. $I_i:=D(e_iA)=\bigoplus_{j\leq 0} D(e_iA_{(-j)})$) is a
projective (resp. an injective) object in $A\mbox{-gr}$. One deduces
that $\{P_i(j)\ |\; 1\leq i \leq n, \ j\in\mathbb{Z}\}$ is a
complete set of pairwise non-isomorphic indecomposable projective
objects in $A\mbox{-gr}$, and $\{I_i(j)\ |\; 1\leq i \leq n, \
j\in\mathbb{Z}\}$ is a complete set of pairwise non-isomorphic
indecomposable injective objects in $A\mbox{-gr}$.
\par
\vskip5pt

Let $0\ne M \in A\mbox{-gr}$. Define $t(M):={\rm max}\{i \in
\mathbb{Z}\ |\; M_{(i)}\neq 0\}$ and $b(M):={\rm min}\{i \in
\mathbb{Z}\ |\; M_{(i)}\neq 0\}$. For a graded $A$-module
$M=\bigoplus_{i \in \mathbb{Z}} M_{(i)}\ne 0$, set ${\rm
top}(M):=M_{(t(M))}$ and ${\rm bot}(M):=M_{(b(M))}$, both of which
are viewed as $A_{(0)}$-modules. Denote by $\Omega^n$ (resp.
$\Omega_{A_{(0)}}$) the $n$-th syzygy functor on $A\mbox{-gr}$
(resp. $A_{(0)}\mbox{-mod}$), \ $n \geq 1$. Similarly we have
$\Omega^{-n}$ and $\Omega_{A_{(0)}}^{-n}$. \par \vskip 5pt

We need the following observation.

\vskip10pt

\begin{lem} Let $M$ be a non-zero non-projective and non-injective graded
$A$-module. With notation above we have \par

 \vskip5pt

$(i)$ Either \ \ $b(\Omega(M))=b(M)$ \ and \ ${\rm
bot}(\Omega(M))=\Omega_{A_{(0)}}({\rm bot}(M))$, \ or \
$b(\Omega(M))> b(M)$.

\vskip5pt

$(i)'$\ \ Either \ $t(\Omega^{-1}(M))=t(M)$\ and \ ${\rm
top}(\Omega^{-1}(M))=\Omega_{A_{(0)}}^{-1}({\rm top}(M))$, \ or\ \
$t(\Omega^{-1}(M)) < t(M)$.
\end{lem}

\noindent{\bf Proof.} \ We only justify $(i)$. Note that ${\rm
rad}(A)={\rm rad}(A_{(0)})\oplus A_{(1)}\oplus \cdots$, and that for
a graded $A$-module $M$, the projective cover $P$ of $M/{{\rm
rad}(A) M}$ in $A\mbox{-mod}$ is graded. It follows that it gives
the projective cover of $M$ in $A\mbox{-gr}$. Since $A$ is
positively-graded, it follows that $b(P)=b(M)$, and that ${\rm
bot}(P)$ is the projective cover of ${\rm bot}(M)$ as
$A_{(0)}$-modules. If ${\rm bot}(P) = {\rm bot}(M)$, then
$b(\Omega(M))> b(M)$. Otherwise, $b(\Omega(M))=b(M)$ \ and \ ${\rm
bot}(\Omega(M))=\Omega_{A_{(0)}}({\rm bot}(M))$.  \hfill
$\blacksquare$

\vskip10pt

\subsection{} Let $A=\bigoplus_{j\geq 0} A_{(j)}$ be a
finite-dimensional positively-graded algebra over $k$. The category
$A\mbox{-gr}$ is said to be \textit{locally representation-finite},
provided that for each $i \in\mathbb{Z}$, the set
 $$\{[M]\; |\; M \mbox{ is
gr-indecomposable such that } M_{(i)}\neq 0\}$$
 is finite, where $[M]$ denote the isoclass in $A\mbox{-gr}$ of the graded module $M$.
By degree-shifts, one sees that $A\mbox{-gr}$ is locally
representation-finite if and only if the set
  $$\{[M]\; |\; M \mbox{ is
gr-indecomposable such that } M_{(0)}\neq 0\}$$ is finite, if and
only if $A\mbox{-gr}$ has only finitely many indecomposable
objects up to degree-shifts.\par \vskip 5pt

If $A$ is in addition self-injective, then $A\mbox{-gr}$ is a
Frobenius category. In fact, we already know that $A\mbox{-gr}$
has enough projective objects and injective objects, and each
indecomposable projective object is of the form $P_i(j)$;  since
$A$ is self-injective, it follows that $P_i$ is injective in
$A\mbox{-mod}$, so is $P_i(j)$ in $A\mbox{-gr}$; similarly, each
$I_i(j)$ is a projective object in $A\mbox{-gr}$.
\par

Note that the stable category $A\underline{\mbox{-gr}}$ is
triangulated (see [Hap], Chap. 1, Sec. 2), with shift functor
induced by $\Omega^{-1}$.

\vskip 10pt

\begin{prop}
 Let $A=\bigoplus_{i\geq 0} A_{(i)}$ be a finite-dimensional positively-graded algebra which is
 self-injective. Assume that
 $ {\rm dim} (A\underline{\mbox{-{\rm gr}}})=0$ and ${\rm gl.dim}(A_{(0)})<
 \infty$. Then $A\mbox{-\rm gr}$ is
 locally representation-finite.
 \end{prop}

 \noindent{\bf Proof.}\quad Since ${\rm dim} (A\underline{\mbox{-{\rm
 gr}}})=0$, it follows that $A\underline{\mbox{-gr}}=\langle X\rangle$ for
 some graded module $X$. Without loss of generality, we may assume that
 $X=\bigoplus_{l=1}^r M^l$, where $M^l$'s are pairwise non-isomorphic
  non-projective gr-indecomposable modules.
  It follows that every
 gr-indecomposable $A$-module is in the set
 $\{ \Omega^i(M^l), \ P_j(i)\ | \ i\in \mathbb{Z}, \ 1\leq l \leq r, \ 1\leq  j\leq n \}$.
Therefore, it suffices to prove that for each $1\leq l \leq
 r$, the set
 $$\{j\in \mathbb{Z}\; | \; \Omega^j(M^l)_{(0)}\neq
 0\}$$
 is finite.\par

For this, assume that ${\rm gl.dim}(A_{(0)})=N$,
 $b(M^l)=j_0$, and $t(M^l)=i_0$. Since ${\rm gl.dim}(A_{(0)}) < \infty$,
 it follows from Lemma 2.1$(i)$ that if $b(\Omega(M))= b(M)$ then
${\rm p.d}({\rm bot}(\Omega(M))) = {\rm p.d}({\rm bot}(M)) - 1$ as
$A_{(0)}$-modules, and otherwise $b(\Omega(M)) > b(M)$. By using
Lemma 2.1$(i)$ repeatedly we have
$$\mbox{ if } \ j\geq {\rm max}\{1, -j_0N\}, \
\mbox{then} \ b(\Omega^j(M^l))>0.$$

Dually, if $j\geq {\rm max}\{1, i_0N\}$,  then
$t(\Omega^{-j}(M^l))< 0$. Note that $b(\Omega^j(M^l))>0$ (resp.
$t(\Omega^{-j}(M^l))<0$) implies that $\Omega^j(M^l)_{(0)} = 0$
(resp. $\Omega^{-j}(M^l)_{(0)}= 0$). It follows that the set
considered above is finite.  \hfill $\blacksquare$

\vskip10pt
\subsection{} Let us recall some related notion
in [BG] and [G]. Let $A$ and $\{e_1, e_2, \cdots, e_n\}$ be the same
as in 2.1, and ${\bf M}$ the full subcategory of $A\mbox{-{\rm gr}}$
consisting of objects $\{P_j(i)\; |\; 1\leq j\leq n,\ i \in
\mathbb{Z}\}$. Then ${\bf M}$ is a locally finite-dimensional in the
sense of [BG]. One may identify $A\mbox{-gr}$ with ${\rm mod}({\bf
M})$ such that a graded $A$-module $M$ is identified with a
contravariant functor sending $P_j(i)$ to $e_jM_{(-i)}$. Now it is
direct to see that $A\mbox{-gr}$ is locally representation-finite if
and only if the category ${\bf M}$ is locally representation-finite
in the sense of [BG], p.337.\par \vskip 5pt

Let us follow [G], p.85-93. Let $G$ be the group $\mathbb{Z}$.
Then $G$ acts freely on ${\bf M}$ by degree-shifts. Moreover,  the
orbit category ${\bf M}/G$ can be identified with the full
subcategory of $A\mbox{-mod}$ consisting of $\{ P_j\; | \; 1\leq
j\leq n\}$. Hence we may identify ${\rm mod}({\bf M}/G)$ with
$A\mbox{-mod}$. With these two identifications, the push-down
functor $F_\lambda: {\rm mod}({\bf M}) \longrightarrow {\rm
mod}({\bf M}/G)$ is nothing but the degree-forgetful functor $U:
A\mbox{-gr} \longrightarrow A\mbox{-mod}$. The following is just a
restatement of Theorem d) in 3.6 of [G].

\begin{lem} Let $k$ be algebraically closed, and $A$
be a finite-dimensional positively-graded $k$-algebra. Assume that
$A\mbox{-{\rm gr}}$ is locally representation-finite. Then the
degree-forgetful functor $U$ is dense, and  hence $A$ is of finite
representation type.
\end{lem}

\vskip10pt
\subsection{Proof of Theorem :}\ If $A$ is an iterated tilted algebra of Dynkin type, then by Theorem
2.10 in [Hap], p.109, we have a triangle-equivalence
$D^b(A\mbox{-mod}) \simeq D^b(kQ\mbox{-mod})$ for some Dynkin quiver
$Q$. Note that $kQ$ is of finite representation type, and that
$D^b(kQ\mbox{-mod}) = \langle M[0]\rangle$, where $M$ is the direct
sum of all the (finitely many) indecomposable $kQ$-modules. It
follows that ${\rm der.dim}(A)={\rm der.dim}(kQ)=0$.

\par

Conversely, if ${\rm dim}D^b(A\mbox{-mod})=0$, it follows from the
fact that $D^b(A\mbox{-mod})$ is Krull-Schmidt (see e.g. Theorem B.2
in Appendix) that $D^b(A\mbox{-mod})$ has only finitely many
indecomposable objects up to shifts. Since $K^b(A\mbox{-proj})$ is a
thick subcategory of $D^b(A\mbox{-mod})$, it follows that
$K^b(A\mbox{-proj})$ has finitely many indecomposable objects up to
shifts. Consequently, ${\rm s.gl.dim}(A) < \infty$ (for the
definition of ${\rm s.gl.dim}(A)$ see B.3 in Appendix).
\par
By Theorem 4.9 in [Hap], p.88, and Lemma 2.4 in [Hap], p.64, we have
an exact embedding
$$F:D^b(A\mbox{-mod})\longrightarrow T(A)\underline{\mbox{-gr}},$$
where $T(A)=A\oplus DA$ is the trivial extension algebra of $A$,
which is graded with ${\rm deg}A=0$ and ${\rm deg}DA=1$. Since ${\rm
gl.dim}A \le {\rm s.gl.dim}(A)-1<\infty$ (see Corollary B.3 in
Appendix), it follows from Theorem 4.9 in [Hap] that the embedding
$F$ is an equivalence. Now by applying Proposition 2.2 to the graded
algebra $T(A)$ we know that $T(A)\mbox{-gr}$ is locally
representation-finite. It follows from Lemma 2.3 that $T(A)$ is of
finite representation type, and then the assertion follows from a
theorem of Assem, Happel, and Rold\`an in [AHR], which says the
trivial extension algebra $T(A)$ is of finite representation type if
and only if $A$ is an iterated tilted algebra of Dynkin type (see
also Theorem 2.1 in [Hap], p.199, and [HW]). \hfill $\blacksquare$
\bibliography{}

\newpage

\centerline {\bf Appendix} \vskip10pt

This appendix includes an exposition on some material we used. They
are well-known, however their proofs seem to be scattered in various
literature.

\vskip10pt\centerline{ A. \ \  Krull-Schmidt categories}

\vskip10pt

This part includes a review of Krull-Schmidt categories.

\par \vskip 5pt

{\bf A.1.} An additive category $\mathcal{C}$ is
\textit{Krull-Schmidt} if any object $X$ has a decomposition $X=X_1
\oplus \cdots \oplus X_n$, such that each $X_i$ is indecomposable
with local endomorphism ring (see [Ri], p.52).

\par \vskip 10pt

Directly by definition, a factor category (see [ARS], p.101) of a
Krull-Schmidt category is Krull-Schmidt.

\vskip10pt

Let $\mathcal{C}$ be an additive category. An idempotent $e =
e^2\in {\rm End}_\mathcal{C}(X)$ \textit{splits}, if there are
morphisms $u: X \longrightarrow Y$ and $v : Y \longrightarrow X$
such that $e=vu$ and ${\rm Id}_Y=uv$. In this case, $u$ (resp.
$v$) is the cokernel (resp. kernel) of ${\rm Id}_X -e$; and ${\rm
End}_\mathcal{C}(Y) \simeq e{\rm End}_\mathcal{C}(X)e$ by sending
$f \in {\rm End}_\mathcal{C}(Y)$ to $vfu$. If in addition ${\rm
Id}_X-e$ splits via $X \stackrel{u'}{\longrightarrow} Y'
\stackrel{v'}{\longrightarrow} X$, then $\binom{u}{u'}: \ X \simeq
Y \oplus Y'$. One can prove directly that an idempotent $e$ splits
if and only if the cokernel of ${\rm Id}_X -e$ exists, if and only
if the kernel of ${\rm Id}_X-e$ exists. It follows that if
$\mathcal{C}$ has cokernels (or kernels) then each idempotent in
$\mathcal{C}$ splits; and that if each idempotent in $\mathcal{C}$
splits, then each idempotent in a full subcategory $\mathcal{D}$
splits if and only if $\mathcal{D}$ is closed under direct
summands.\par \vskip 5pt

\vskip10pt

A ring $R$ is \textit{semiperfect} if $R/{\rm rad}(R)$ is
semisimple and any idempotent in $R/{\rm rad}(R)$ can be lifted to
$R$, where ${\rm rad}(R)$ is the Jacobson radical. \vskip10pt

\noindent {\bf Theorem A.1.}\  \emph{An additive category
$\mathcal C$ is Krull-Schmidt if and only if any idempotent in
$\mathcal C$ splits, and ${\rm End}_\mathcal{C}(X)$ is semiperfect
for any $X \in \mathcal{C}$.}

\emph{In this case, any object has a unique (up to order) direct
decomposition into indecomposables.}

\vskip5pt

\noindent{\bf Proof.}\quad For $X \in \mathcal{C}$ denote by ${\rm
add}X$ the full subcategory of the direct summands of finite direct
sums of copies of $X$, and set $R:={\rm End}_\mathcal{C}(X)^{op}$.
Let $R{\rm \mbox{-}proj}$ \ denote the category of
finitely-generated projective left $R$-modules. Consider the
fully-faithful functor
\begin{align*}
\Phi_X:= {\rm Hom}_\mathcal{C}(X, -): \ {\rm add}X \longrightarrow
R\mbox{-proj}.
\end{align*}

Assume that $\mathcal{C}$ is Krull-Schmidt. Then $X=X_1 \oplus
\cdots \oplus X_n$ with each $X_i$ indecomposable and ${\rm
End}_\mathcal{C}(X_i)$ local. Set $P_i:=\Phi_X(X_i)$. Then $_R
R=P_1 \oplus \cdots \oplus P_n$ with ${\rm End}_R(P_i) \simeq {\rm
End}_\mathcal{C}(X_i)$ local. Thus $R$ is semiperfect by Theorem
27.6(b) in [AF], and so is ${\rm End}_\mathcal{C}(X)=R^{op}$. Note
that every object $P \in R\mbox{-proj}$ is a direct sum of
finitely many $P_i$'s: in fact, note that $\{S_i:=P_i/{\rm
rad}(P_i)\}_{1 \leq i \leq n}$ is the set of pairwise
non-isomorphic simple $R$-modules and that the projection $P
\longrightarrow P/{\rm rad}(P)= \bigoplus\limits_i S_i^{m_i}$ is a
projective cover, thus $P \simeq \bigoplus\limits_i P_i^{m_i}$. It
follows that $P$ is essentially contained in the image of
$\Phi_X$, and hence $\Phi_X$ is an equivalence. Consider
$R\mbox{-Mod}$, the category of left $R$-modules. Since
$R\mbox{-Mod}$ is abelian, it follows that any idempotent in
$R\mbox{-Mod}$ splits. Since $R\mbox{-proj}$ is a full subcategory
of $R\mbox{-Mod}$ closed under direct summands, it follows that
any idempotent in $R\mbox{-proj}$ splits. So each idempotent in
${\rm add}(X)$ splits. This proves that any idempotent in
$\mathcal C$ splits.
\par

Conversely, assume that each idempotent in $\mathcal{C}$ splits
and $R^{op}={\rm End}_\mathcal{C}(X)$ is semiperfect for each $X$.
Then again by Theorem 27.6(b) in [AF] we have $R=Re_1 \oplus
\cdots \oplus Re_n$ where each $e_i$ is idempotent such that
$e_iRe_i$ is local. Since $1= e_1 + \cdots + e_n$ and $e_i$ splits
in $\mathcal{C}$ via, say $X \stackrel{u_i}{\longrightarrow} Y_i
\stackrel{v_i}{\longrightarrow} X$, it follows that $X \simeq
Y_1\oplus \cdots \oplus Y_n$ via the morphism $(u_1, \cdots,
u_n)^t$ with inverse $(v_1, \cdots, v_n)$. Note that ${\rm
End}_\mathcal{C}(Y_i) \simeq e_i {\rm
End}_\mathcal{C}(X)e_i=(e_iRe_i)^{op}$ is local. This proves that
$\mathcal{C}$ is Krull-Schmidt. \par

For the last statement, it suffices to show the uniqueness of
decomposition in ${\rm add}X$ for each $X$. This follows from the
fact that $\Phi_X$ is an equivalence, since the uniqueness of
decomposition in $R\mbox{-proj}$ is well known by Azumaya's
theorem (see e.g. Theorem 12.6(2) in [AF]). This completes the
proof. \hfill $\blacksquare$

\vskip10pt

{\bf A.2.} Let $k$ be a field. An additive category $\mathcal{C}$
is a Hom-finite $k$-category if ${\rm Hom}_\mathcal{C}(X, Y)$ is
finite-dimensional $k$-space for any $X, Y \in \mathcal{C}$, or
equivalently, ${\rm End}_\mathcal{C}(X)$ is a finite-dimensional
$k$-algebra for any object $X$.

\vskip10pt

\noindent {\bf Corollary A.2.}\ \emph{Let $\mathcal{C}$ be a
Hom-finite $k$-category. Then the following are equivalent.}

\vskip5pt

 $(i)$\ \  $\mathcal{C}$ \emph{is Krull-Schmidt.}\par
 \vskip5pt
 $(ii)$ \ \emph{Each idempotent in $\mathcal{C}$ splits.}\par
 \vskip5pt
 $(iii)$\ \emph{For any indecomposable $X \in
\mathcal{C}$, ${\rm End}_{\mathcal{C}}(X)$ has no non-trivial
idempotents.}\par

 \vskip 10pt

\noindent {\bf Remark A.2.} \ \emph{By Corollary A.2 $(ii)$, a
Hom-finite abelian $k$-category is Krull-Schmidt.}

\emph{In particular, the category of coherent sheaves on a
complete variety is Krull-Schmidt (see [At], Theorem 2(i)).}

\vskip20pt

\centerline {B. \ \ Homotopically-minimal complexes}

\vskip10pt

In this part $A$ is a finite-dimensional algebra over a field $k$.

\vskip10pt

{\bf B.1.} \ A complex $P^\bullet = (P^n, d^n)\in
C(A\mbox{-proj})$ is called \textit{homotopically-minimal}
provided that a chain map $\phi^\bullet: P^\bullet \longrightarrow
P^\bullet$ is an isomorphism if and only if it is an isomorphism
in $K(A\mbox{-proj})$ (see [K]).\par \vskip 5pt

Applying Lemma B.1 and Proposition B.2 in [K], and duality,  we have

\vskip10pt

\noindent {\bf Proposition B.1.} {\rm  (Krause)} \emph{Let
$P^\bullet=(P^n, d^n) \in C(A\mbox{-{\rm proj}})$. The following
statements are equivalent.}
\par
\vskip5pt $(i)$\quad \emph{The complex $P^\bullet$ is
homotopically-minimal.}\par

 \vskip5pt

 $(ii)$ \quad \emph{Each differential $d^n$ factors through ${\rm
rad}(P^{n+1})$}. \par

 \vskip5pt

 $(iii)$\quad \emph{The complex $P^\bullet$ has no non-zero direct
summands in $C(A\mbox{-{\rm proj}})$ which are null-homotopic.}
\par \vskip5pt

\emph{Moreover, in $C(A\mbox{-{\rm proj}})$ every complex
$P^\bullet$ has a decomposition $P^\bullet=P'^\bullet \oplus
P''^\bullet$ such that $P'^\bullet$ is homotopically-minimal and
$P''^\bullet$ is null-homotopic.}

\vskip10pt

{\bf B.2.} \ For $P^\bullet \in C(A\mbox{-proj})$, consider the
ideal of ${\rm End}_{C(A\mbox{-proj})}(P^\bullet)$:
\begin{align*}{\rm
Htp}(P^\bullet)=\{\phi^\bullet: P^\bullet\longrightarrow
P^\bullet\ |\; \phi^\bullet \mbox{ is homotopic to zero} \}.
\end{align*}
\noindent {\bf Lemma B.2.}\  \emph{Assume ${\rm rad}^l (A)=0$. Let
$P^\bullet$ be homotopically-minimal. Then ${\rm
(Htp}(P^\bullet))^l=0$.}

\vskip5pt

\noindent{\bf Proof.}\quad Let $\phi^\bullet \in {\rm
Htp}(P^\bullet)$ with homotopy $\{h^n\}$. Then $\phi^n=d^{n-1}h^n+
h^{n+1}d^{n}$. Since by assumption both $d^{n-1}$ and $d^{n}$ factor
through radicals, it follows that $\phi^n$ factors through ${\rm
rad}P^n$. Therefore, for $k\geq 1$ morphisms in $({\rm
Htp}(P^\bullet))^k$ factor through the $k$-th radicals. So the
assertion follows from ${\rm rad}^l(A)=0$.\hfill $\blacksquare$

\vskip10pt

Denote by $C^{-, b}(A\mbox{-proj})$ the category of bounded above
complexes of projective modules with finitely many non-zero
cohomologies, and by $K^{-,b}(A\mbox{-proj})$ its homotopy
category. It is well-known that there is a triangle-equivalence
${\bf p}: D^b(A\mbox{-mod}) \simeq K^{-,b}(A\mbox{-proj})$.

\vskip10pt

The following result is can be deduced from Corollary 2.10 in [BS].
See also [BD]. \vskip10pt

\noindent {\bf Theorem B.2.} \  \emph{The bounded derived category
$D^b(A\mbox{-mod})$ is Krull-Schmidt.} \vskip5pt

\noindent{\bf Proof.}\quad Clearly $D^b(A\mbox{-mod})$ is
Hom-finite. By Corollary A.2 it suffices to show that ${\rm
End}_{D^b(A\mbox{-mod})}(X^\bullet)$ has no non-trivial idempotents,
for any indecomposable  $X^\bullet$.\par

By Proposition B.1 we may assume that $P^\bullet: = {\bf
p}X^\bullet$ is homotopically-minimal. Since $P^\bullet$ is
indecomposable in $K^{-,b}(A\mbox{-proj})$, it follows from
Proposition B.1$(iii)$ that $P^\bullet$ is indecomposable in
$C(A\mbox{-proj})$. Since idempotents in $C(A\mbox{-proj})$ split,
it follows that ${\rm End}_{C(A\mbox{-proj})}(P^\bullet)$ has no
non-trivial idempotents. Note that
\begin{align*}
{\rm End}_{D^b(A\mbox{-mod})}(X^\bullet)={\rm
End}_{K^{-,b}(A\mbox{-proj})}(P^\bullet)={\rm
End}_{C(A\mbox{-proj})}(P^\bullet)/{{\rm Htp}(P^\bullet)}.
\end{align*}
Since by Lemma B.2 ${\rm Htp}(P^\bullet)$ is a nilpotent ideal, it
follows that any idempotent in the quotient algebra ${\rm
End}_{C(A\mbox{-proj})}(P^\bullet)/{{\rm Htp}(P^\bullet)}$ lifts to
${\rm End}_{C(A\mbox{-proj})}(P^\bullet)$. Therefore ${\rm
End}_{C(A\mbox{-proj})}(P^\bullet)/{{\rm Htp}(P^\bullet)}$ has no
non-trivial idempotents. \hfill $\blacksquare$

\vskip10pt

{\bf B.3.} \ For $X^\bullet=(X^n, d^n)$ in $C^b(A\mbox{-mod})$,
define the \textit{width} $w(X^\bullet)$ of $X^\bullet$ to be the
cardinality of $\{n \in \mathbb{Z}| \; X^n \neq 0\}$. The strong
global dimension ${\rm s.gl.dim}(A)$ of $A$ is defined by (see
[S])
\begin{align*}
{\rm s.gl.dim}(A):= {\rm sup}\{w(X^\bullet)\ | \; X^\bullet \mbox{
is indecomposable in } C^b(A\mbox{-proj})\}.
\end{align*}

By Proposition B.1 an indecomposable $X^\bullet$ in
$C^b(A\mbox{-proj})$ is either homotopically-minimal, or
null-homotopic (thus it is of the form $ \cdots \longrightarrow 0
\longrightarrow P \stackrel{{\rm Id}}{\longrightarrow} P
\longrightarrow  0 \longrightarrow \cdots$, for some indecomposable
projective $A$-module $P$). So we have
\begin{align*}
{\rm s.gl.dim}(A)= {\rm sup}\{2, \ w(P^\bullet)&\ | \; P^\bullet
\mbox{  is homotopically-minimal} \\
&\mbox{and indecomposable in}\  C^b(A\mbox{-proj})\}.
\end{align*}

Let $M$ be an indecomposable $A$-module with minimal projective
resolution $P^\bullet\stackrel{\varepsilon}{\longrightarrow} M $.
Denote by $\tau^{\geq -m}P^\bullet$ the brutal truncation of
$P^\bullet$, $m\geq 1$. By Proposition B.1$(ii)$ $\tau^{\geq
-m}P^\bullet$ is homotopically-minimal.

If $\tau^{\geq -m}P^\bullet = P'^\bullet\oplus Q^\bullet$ in
$C^b(A\mbox{-proj})$ with $P^\bullet=(P^n, d^n)$, $P'^\bullet=(P'^n,
\delta^n)$, and $Q^\bullet=(Q^n, \partial^n)$, then both
$P'^\bullet$ and $Q^\bullet$ are homotopically-minimal. Assume that
$P'^0\neq 0$, and set $t_0: ={\rm max}\{ t \in \mathbb{Z}\ |\;
Q^t\neq 0\}$. Then $-m\leq t_0\leq 0$. Since $M$ is indecomposable
and both $P'^\bullet$ and $Q^\bullet$ are homotopically-minimal, it
follows that $t_0\ne 0$, and hence $Q^{t_0} \subseteq \mbox{{\rm
Ker}} d^{t_0} \subseteq \mbox{{\rm rad}}(P^{t_0}) = \mbox{{\rm
rad}}(P'^{t_0}\oplus Q^{t_0})$, a contradiction. This proves

\vskip10pt

\noindent {\bf Lemma B.3.}\  \emph{The complex $\tau^{\geq
-m}P^\bullet$ is homotopically-minimal and indecomposable in
$C^b(A\mbox{-{\rm proj}})$.} \vskip5pt

\par \vskip 10pt

As a consequence we have

\vskip10pt

\noindent {\bf Corollary B.3.} \ ([S], p.541) \emph{Let $A$ be a
finite-dimensional algebra. Then
\begin{align*}
{\rm s.gl.dim}(A) \geq {\rm max}(2, 1+ {\rm gl.dim}(A)).
\end{align*}}

\end{document}